\documentclass[12pt]{amsart}
\usepackage{amssymb}
\usepackage{amscd}

\newcommand{\mm}{\mathfrak m}

\newcommand{\Z}{\mathbb{Z}}
\newcommand{\R}{\mathbb{R}}
\newcommand{\N}{\mathbb{N}}

\newcommand{\Mcc}{m}
\newcommand{\Ncc}{n}
\newcommand{\Tc}{\mathcal{T}}

\DeclareMathOperator{\stand}{\mathcal S}

\DeclareMathOperator{\ini}{in}
\DeclareMathOperator{\Ker}{Ker}

\DeclareMathOperator{\rank}{rank}

\DeclareMathOperator{\supp}{supp}

\DeclareMathOperator{\id}{id}

\DeclareMathOperator{\tensor}{\otimes}

\DeclareMathOperator{\pnt}{\raise 0.5mm \hbox{\large\bf.}}

\DeclareMathOperator{\str}{{\rm star}}
\DeclareMathOperator{\rk}{{\rm rank}}
\DeclareMathOperator{\depth}{{\rm depth}}

\DeclareMathOperator{\op}{op}

\DeclareMathOperator{\cone}{{\rm cone}}

\newtheorem{thm}{\bf Theorem}[section]
\newtheorem{lem}[thm]{\bf Lemma}
\newtheorem{cor}[thm]{\bf Corollary}

\theoremstyle{definition}
\newtheorem{defn}[thm]{\bf Definition}

\newtheorem{ex}[thm]{\bf Example}

\title{On Algebras associated to partially ordered sets}

\author{Morten Brun}
\address{FB Mathematik/Informatik, Universit\"at Osnabr\"uck, 49069 Osnabr\"uck, Germany}
\email{brun@mathematik.uni-osnabrueck.de}

\author{Tim R\"omer}
\address{FB Mathematik/Informatik, Universit\"at Osnabr\"uck, 49069 Osnabr\"uck, Germany}
\email{troemer@mathematik.uni-osnabrueck.de}

\begin{document}

\begin{abstract}
We continue the study \cite{BBR05}
on sheaves of rings on finite posets.
We present examples where the ring of global sections
coincide with toric faces rings, quotients of a polynomial ring
by a monomial ideal and algebras with straightening laws.
We prove a rank-selection theorem
which generalizes the well-known
rank-selection theorem of Stanley--Reisner rings.
Finally, we determine an explicit
presentation of certain global rings of sections.
\end{abstract}
\maketitle
%
%
%

\section{Introduction}
In the present paper we continue the study \cite{BBR05} on
sheaves of commutative rings on finite partially ordered sets (posets for short).
Recall that a finite poset $P$ can be considered as a topological
space
with the {\em Alexandrov topology} \cite{Al56}, that is, the
topology where the open sets are the
subsets
$U$ of $P$ such that $y \in U$ and $x \le y$ implies $x \in U$.
Let  $R$ be a $\Z^m$-graded commutative Noetherian ring.
A sheaf $\Tc$ of $\Z^m$-graded $R$-algebras on $P$
is uniquely determined by its stalks $(T_x)_{x \in P}$ and the
homomorphisms $T_{xy} \colon T_y \to T_x$ for $x \le y$
in $P$ induced by the restriction maps of $\Tc$. More precisely, the
stalks of $\Tc$ form an
{\em $RP$-algebra} $T$, that is, a system $(T_x)_{x \in P}$ of
$\Z^m$-graded $R$-algebras and
$\Z^m$-graded homomorphisms $T_{xy} \colon T_y \to T_x$ for $x \le y$
in $P$ with
the property that $T_{xx}$ is the identity on $T_x$ and that $T_{xy}
\circ T_{yz} = T_{xz}$ for every $x \le y \le z$ in $P$. Moreover,
every $RP$-algebra $T$ corresponds to a sheaf $\Tc$.
We study the (inverse) {\em limit} $\lim_P T \cong H^0(P,\Tc)$ of $T$
and call it the {\em ring of global sections of $T$}.
Under this name it was for example studied by Yuzvinsky \cite{YU87}
and Caijun \cite{CA97}.
The $RP$-algebra $T$ is called {\em flasque} if $\Tc$ is flasque, that
is,
if the restriction $\lim_U T \to \lim_V T$
is surjective for every inclusion
$V \subseteq U$ of open subsets of $P$.
The case where $R = K$ is a field is of particular interest.

The paper is organized as follows: In Section \ref{prereq} we explain
the
terminology used throughout the paper.
In \cite{BBR05} it was shown that rings of global sections of
$KP$-algebras appear
naturally in algebraic combinatoric
and commutative algebra. In Section \ref{examples}
we
present classes of $KP$-algebras
associated to toric face rings (as introduced in \cite{BRRO04}),
monomial ideals and algebras with straightening law (ASL's for
short).

In Section \ref{aslsheaf} we give a criterion on
a $KP$-algebra with stalks given by ASL's
to have an ASL as ring of global sections.
As an application we generalize
a construction of Stanley \cite{ST91} who defined
the face ring of a simplicial poset.
Here we consider more generally a locally distributive lattice,
i.e.
$P$ has a terminal element $0_P$ and
for all $z \in P$ the interval
$[0_P,z]$ is a distributive lattice.
We construct
an $KP$-algebra $T$ such that
$\lim T \cong K[P]/I_{T}$
where $K[P]$ is a polynomial ring on the elements $x \in P$
and
the ideal
$I_{T}$
is generated
by
\begin{eqnarray*}
x y &-& x\wedge y \sum_{z } z
\end{eqnarray*}
where $z$ ranges over all minimal upper bounds for $x,y$ in $P$.
The sum is defined to be zero if there are no such elements.

After the variety of examples we
prove in Section \ref{quotients}
a rank selection theorem
in the tradition of
Duval \cite{DU97},
Hibi \cite{HI91}, Munkres \cite{MU84}
and Stanley \cite{ST79}.
As an application of one of our main results
Theorem \ref{remCMthm}
we consider the following situation.
We define
the
{\em $i$th skeleton} of a finite poset $P$ to be
the sub-poset
$P^{(i)}=\{y \in P:\rk(y) \leq i\}.$
If for example
$P=P(\Sigma)$ is the face poset of a rational pointed fan $\Sigma$
in $\R^m$
and $T=T(\Sigma)$ is the $KP$-algebra as constructed in
Example \ref{exfan}, then
we obtain that
$$
\depth \lim_P T
=
\max\{i \in \N : i \le \rank(P) \text{ and } \lim_{P^{(i)}} T \text{
  is Cohen--Macaulay} \}.
$$
If $P$ is the face poset of a simplicial complex $\Delta$ and $T$ is
the Stanley--Reisner $KP$-algebra associated to $\Delta$, then this
is a well-known result, because
$\lim T\cong K[\Delta]$
is the Stanley--Reisner ring associated to $\Delta$
and $\lim_{P^{(i)}} T\cong K[\Delta^{(i)}]$
is the Stanley--Reisner ring associated to the $i$th skeleton
$\Delta^{(i)}$ of $\Delta$ in this situation.

Finally, in Section \ref{sectioned} we study presentations of $KP$-algebras.
Many such $KP$-algebras
have an extra piece of
structure: there exist compatible splittings of the restriction maps
considered as homomorphisms of abelian groups. These splitting enable
us to present the associated rings of global sections of
$T$ by an ideal $I$ of the form
$I = I_{\Tc} + I_{\Delta}$
in a polynomial ring,
where $I_{\Delta}$ is the Stanley--Reisner ideal of a simplicial
complex $\Delta$ related to the order complex of the underlying poset,
and $I_{\Tc}$ is the sum of the defining ideals for the stalks $T_x$
of
$T$. Given such a presentation, the main result of \cite{BBR05}
implies that the cohomology of the simplicial
complex $\Delta$ carries enough information to determine the depth and
dimension of the ring of global sections $\lim_P T$ of $T$.

\section{Prerequisites}
\label{prereq}
In this paper $P = (P,\le)$
denotes a {\em partially ordered set} (poset for short).
We write $x < y$ if $x \neq y$ and $x \le y$.
The opposite poset $P^{\op}$ has the same
underlying set as
$P$ and the opposite partial order.
If $P$ contains a unique maximal element $1_P$, then it
is called the {\em initial element} of $P$.
Analogously a unique minimal element $0_P$ of $P$ is called a {\em terminal
  element} of $P$.
The poset $\widehat P$ is obtained from $P$ by adjoining disjoint
initial- and terminal elements $1_{\widehat P}$ and $0_{\widehat P}$.
The closed interval $[x,y]$ is the set of elements between $x$ and $y$
in $P$, that is, elements $z \in P$ with $x \le z \le y$.
We also consider the half-open interval $[x,y) = [x,y]
  \setminus\{y\}$, the half-open interval $(x,y] = [x,y]
\setminus\{x\}$ and the open interval $(x,y) = [x,y] \setminus\{x,y\}$.
A {\em chain} is a
totally ordered subset of $P$.
The rank of $P$, denoted $\rank(P)$,
is the supremum of the numbers $|C|-1$ for $C$ a finite chain in
$P$.
The rank function on $P$ is defined by $\rank(x) = \rank((0_{\widehat
  P},x])$ for $x \in P$.
The {\em face poset} of a simplicial complex $\Delta$ is the poset
$P(\Delta) =
(\Delta, \subseteq)$ of {\em faces} in $\Delta$ ordered by inclusion.
The {\em order complex}
$\Delta(P)$ of $P$ is the simplicial complex
consisting of the chains in $P$ ordered by inclusion.
(For more details
see
Bruns--Herzog \cite{BRHE98} and Stanley \cite{ST96,ST99}.)

Let  $R$ be a $\Z^m$-graded commutative Noetherian ring and let $T$ be
an $RP$-algebra.
If the rings $T_x$ have property $E$ and the homomorphisms $T_{xy}$
are compatible with property $E$
we say that
$T$ is an $RP$-algebra of algebras with property $E$.
For example $T$ is an $RP$-algebra of Cohen--Macaulay rings,
if the rings $T_x$ are Cohen--Macaulay.
If there exists exactly one
$\Z^m$-graded ideal $\mm$ in $R$ such that $R/\mm$ is a field, then we
call
$R$ a {\em $\Z^m$-graded local ring}.
Note that $R/\mm$ is concentrated in degree zero, and thus $\mm$
contains every element of $R$ of non-zero degree.
In particular, if $K = R_0$ is a field, then
$\mm = \bigoplus_{a \in \Z^m \setminus \{0\}} R_a$.
We denote by $H^i_\mm(N)$ the local cohomology modules of
a finitely generated $\Z^m$-graded $R$-module $N$ with respect
to $\mm$ (see \cite{BRHE98} for details). Since $\mm$ is a maximal
ideal the definition of graded local cohomology modules and usual
local cohomology modules coincide. Note that $H^i_\mm(N)$ is a
$\Z^m$-graded $R$-module.
\section{Examples of $RP$-algebras}
\label{examples}
In this section we present $RP$-algebras appearing
naturally in algebraic combinatoric and commutative algebra.
Of course every $K$-algebra $S$ equals $\lim_P T$ for
the poset $P=\{x\}$ consisting of one element
and the $KP$-algebra $T$ with $T_x=S$.
But this gives no new information about $S$.
The goal is to choose
a finite poset $P$ and a
suitable $KP$-algebra $T$
such that the stalks $T_x$ are as nice as possible (e.g. Cohen--Macaulay rings)
and to characterize ring properties of $\lim_P T$
in terms of combinatorial properties of $P$
and ring properties of the rings $T_x$ for $x \in P$.
Our first example goes back to a construction of Stanley \cite{ST87}.

\begin{ex}[Toric face rings]
\label{exfan}
We consider a rational fan $\Sigma$ in $\R^m$,
that is, $\Sigma$ is a finite collection of rational cones in $\R^m$
such that for $C' \subseteq C$ with $C\in \Sigma$
we have that $C'$ is a face of $C$ if and only if $C' \in \Sigma$,
and if $C,C' \in \Sigma$, then $C\cap C'$ is a common face of $C$ and $C'$.
For each cone $C \in \Sigma$
we choose an affine monoid $M_C \subseteq \Z^m$
such that:
\begin{enumerate}
\item
$\cone(M_C)=C$ for $C \in \Sigma$;
\item
If $C,C' \in \Sigma$, $C' \subseteq C$, then $M_{C'} = M_C \cap C'$.
\end{enumerate}
Observe that we do not require that $M_C$ is normal.
Let $P(\Sigma) = (\Sigma, \subseteq)$ be the partially ordered set
of faces of $\Sigma$ ordered by inclusion and
let $K$ be a field.
For $C \in P(\Sigma)$ we let $T_{C}$ denote the affine monoid ring $K[M_C]$.
For $C' \subseteq C$ in $P(\Sigma)$ the homomorphisms $T_{C'C}\colon
T_{C} \to T_{C'}$
are induced by the natural face projections $K[C \cap
  \Z^m] \to K[C' \cap \Z^m]$.
This is a $\Z^m$-graded $KP$-algebra, and
with $R=\lim T$ it is an $RP$-algebra of
$\Z^m$-graded $R$-algebras.
We write $T = T(\Sigma)$ if $M_C=C \cap \Z^m$, for every $C \in
\Sigma$.
In this case the ring $\lim T(\Sigma)$ is called the {\em toric face
  ring} of $\Sigma$, and it
was studied by Stanley \cite{ST87}, Bruns--Gubeladze \cite{BRGU01}
and Brun--R\"omer \cite{BRRO04}.
\end{ex}

Stanley--Reisner rings of finite simplicial complexes are
covered by the above example: To a simplicial complex $\Delta$ we
associate a fan $\Sigma$ with the same face poset as $\Delta$, and we
can choose
$M_C$ isomorphic to $\N^{\dim(C)}$ for every
$C \in \Sigma$.
In this situation
we write $P(\Delta) = P(\Sigma)$ and call $T(\Delta) = T(\Sigma)$ the
{\em Stanley--Reisner} $KP$-algebra.
The ring $\lim T(\Delta)$ is isomorphic to the
Stanley--Reisner ring $K[\Delta]$. (See \cite[Example 5.2]{BBR05}.)

Our next example shows that rings defined by monomial ideals
give rise to $KP$-algebras.
\begin{ex}[Monomial ideals]
\label{exmulti}
Let $R=K[x_1,\dots,x_n]$ be the polynomial ring over a field $K$
with the usual $\Z^n$-grading.
Recall that irreducible monomial ideals in $S$ are of the form
$\mm^{b}=(x_1^{b_1}, \dots, x_n^{b_n})$ for $0 \ne b \in \N^n$.
A monomial ideal $I$ has a unique irredundant
irreducible decomposition of the form
$$
I=\mm^{b^1} \cap \cdots \cap \mm^{b^t} \text{ for } b^1,\dots, b^t \in \N^n.
$$
(See \cite[Section 5.2]{MIST} for details.)
For each subset $C$ of $[t] = \{1,\dots,t\}$ there exists a unique maximal
subset $\overline C \subseteq [t]$
with the property that
${\sum_{i \in C}
  \mm^{b^i}} = {\sum_{i \in \overline C}
  \mm^{b^i}}$.
Let
$P(I)$ denote the poset of the subsets $\overline C$ of $[t]$ for
$\emptyset \ne C
\subseteq [t]$
ordered by
reverse inclusion.
For $C \in P(I)$ let $T(I)_{C}=S/{\sum_{i \in C} \mm^{b^i}}$
and for $C, D \in P(I)$ with $C \subseteq D$
we define
$T(I)_{D C}: T_{C} \to T_{D}$ to be the natural projection map.
Now:
\begin{enumerate}
\item
$T(I)$ is a $\Z^n$-graded $RP(I)$-algebra.
\item
By using Lemma 7.3 of \cite{BBR05} we see that
$T(I)$ is flasque.
\item
$T(I)_{C}$ is a complete intersection and thus Cohen--Macaulay for
$C \in P(I)$.
\item
Since we have chosen an irredundant irreducible decomposition, the
one-point sets $\{i\}$ are the maximal elements of $P(I)$, and so
$\lim
T(I) \cong S/I$.
\end{enumerate}
\end{ex}

The result \cite[Theorem 6.2]{BBR05}
has
the following corollary.
(See \cite[Theorem 1]{TA05} for a more general description of
$H^i_\mm(S/I)$.)
\begin{cor}
Let $I = \mm^{b^1} \cap \cdots \cap \mm^{b^t} \subset
K[x_1,\dots,x_n]$ be an irredundant irreducible decomposition of a
monomial ideal and let $d_{C}=\dim T(I)_{C}$. If
$d_{C} > d_{D}$ for $C,D \in P(I)$ with $D \subset C$ then
$H^i_\mm(K[x_1,\dots,x_n]/I)$ is isomorphic to
$$\bigoplus_{F \in P(I)}
\widetilde H^{i - d_{F} -1}((F,1_{\widehat P(I)});K)
\tensor_K
H^{d_{F}}_\mm(K[x_1,\dots,x_n]/{\sum_{i \in F} m^{b^i}})
$$
as a $\Z^n$-graded $K$-vector space.
\end{cor}

For the next example we
recall
the notion of an {\em algebra with straightening law} (ASL for short).
We call a function $\Mcc \colon P \to \N$ (i.e. $m \in \N^P$)
a {\em monomial} on $P$, and consider it as an element of the
polynomial ring $K[P]$ on the elements of $P$.
The {\em support} of $\Mcc$
is the set $\supp(\Mcc)=\{x \in P: \Mcc(x)\neq 0 \}$.
The monomial is called a {\em standard monomial}
if $\supp(\Mcc)$ is a chain in $P$.
Let $R$ be an algebra over a field $K$ with an
injection $\phi\colon P \to R$.
We associate to each monomial $\Mcc$
the element
$
\phi(\Mcc)=\prod_{x \in P} \phi(x)^{\Mcc(x)} \in R$
and call $\phi(\Mcc)$ a {\em monomial in $R$}.
Following \cite{COEICL} we call $R$ an {\em algebra with straightening
  law} (ASL for short)
on $P^{\op}$ over $K$ if
\begin{enumerate}
\item[ASL1]
The set of standard monomials is a $K$-basis of $R$, and
\item[ASL2]
if $x$ and $y$ are incomparable in $P$
and if
\begin{equation}
\label{asleqn}
xy= \sum_{\Mcc} r_{\Mcc,xy} \phi(\Mcc)
\end{equation}
is the unique representation of $xy$ as a linear combination of
standard monomials
guaranteed by ASL1,
then $r_{\Mcc,xy} \ne 0$ implies that
the maximal element $x_{\Mcc}$ of $\supp(\Mcc)$
satisfies $x < x_{\Mcc}$ and $y<x_{\Mcc}$.
\end{enumerate}
Note that for technical reasons we work with $P^{\op}$ instead of $P$.
By ASL1 every monomial $\phi(\Ncc)$ in $R$ has a unique representation
of the form $\phi(\Ncc) = \sum_{\Mcc \in \N^P} r_{\Mcc,\Ncc}
\phi(\Mcc)$, where $r_{\Mcc,\Ncc} \ne 0$ implies that $\Mcc$ is a
standard monomial.
We write
$$r(\Ncc)= \sum_{\Mcc \in \N^P} r_{\Mcc,\Ncc} \Mcc$$
for the element in the
polynomial ring
$K[P]$ associated
to
this representation of $\phi(\Mcc)$.
\begin{ex}[ASL]
\label{asl}
Let $R$ be an ASL on  $P^{\op}$.
Then
$R=K[P]/I_P$ where
the ideal $I_P \subseteq K[P]$ is generated by the straightening relations
$\Ncc - r(\Ncc)$.
For $Q' \subseteq Q$ an inclusion of subsets of $P$
we
consider the projection map
$$
p_{Q'Q} \colon K[Q] \to K[Q'],\
x \mapsto
\begin{cases}
x & \text{if } x \in Q',\\
0 & \text{if } x \not\in Q'.
\end{cases}
$$
We let $I_Q = p_{QP}(I_P) \subseteq K[Q]$, we define $T_Q = K[Q]/I_Q$
and we let $T_{Q'Q} \colon T_Q \to T_{Q'}$ denote the map induced by
$p_{Q'Q}$. For $Q$ open in $P$ the ideal $I_Q$ is generated by the
elements $\Ncc - p_{QP} (r(\Ncc))$ for $\Ncc \in \N^Q$.
By \cite[Proposition 1.2 (b)]{COEICL} (with $I = P \setminus Q$)
the $K$-algebra
$T_Q$ is an ASL on $Q^{\op}$.

Now we let
$T_x=T_{(0_{\widehat P},x]}$ for $x\in P$
and $T_{xy}=T_{(0_{\widehat P},y](0_{\widehat P},x]}$
for $x,y \in P$ with $x\leq y$.
It follows that $T$ is an $RP$-algebra.
Using that $T_Q$ is an ASL on $Q^{\op}$ if $Q \subseteq P$ is open,
we see by working with standard monomials
that $\lim_P T \cong T_P$.
\end{ex}


\section{$KP$-algebras with an ASL-structure}
\label{aslsheaf}
We saw in Example \ref{asl}
that an ASL $R$ can be seen as the ring of global sections of
a suitable $KP$-algebra.
Sometimes it is possible to reverse this construction
in the sense that giving locally ASL's one can construct
a $KP$-algebra such that the global ring of sections
is again an ASL with prescribed presentation.

Recall that $K[P] = K[x : x \in P]$ denotes the
polynomial algebra on the elements $x$ of $P$.
At first we present a criterion which ensures that a given
$KP$-algebra has an ASL-structure.

\begin{thm}
\label{aslconstr}
Let $T$ be a $KP$-algebra on a finite poset $P$ such
that $T_x$ is an ASL on $(0_{\widehat P},x]^{\op}$ for every
$x \in P$.
Let $\phi_x\colon K[(0_{\widehat P},x]] \to T_x$ be the natural
projection induced by the ASL-structure on $T_x$ for $x \in P$.
If the diagram
$$
\begin{CD}
    K[(0_{\widehat P},y]] @>{\phi_y}>> T_y \\
        @V{p_{(0_{\widehat P},x](0_{\widehat P},y]}}VV @V{T_{xy}}VV \\
    K[(0_{\widehat P},x]] @>{\phi_x}>> T_x
\end{CD}
$$
commutes for every $x,y \in P$ with $x \le y$, then $\lim_P T$ is an ASL on
$P^{\op}$.

In particular, if $\stand \subseteq \N^P$ denotes the set of standard monomials in
$K[P]$ then
$\lim_P T \cong K[P]/I_P$,
where the ideal $I_P$ is generated by relations
of the form
$$
\Ncc - \sum_{\Mcc \in \stand} r_{\Mcc, \Ncc} \Mcc
\text{ for } \Ncc \in \N^P.
$$
\end{thm}
\begin{proof}
Let $\phi \colon K[P] \to \lim_P T$ denote the homomorphism
induced by the natural projections
$K[P] \to K[(0_{\widehat P},x]]$ for $x \in P$.

We claim that $\phi$ is surjective and moreover $\lim_P T$
is generated as a $K$-vector space by the elements
$\phi(\Mcc)$ with $\Mcc \in \stand$.

Let $s=(s_x)_{x \in P} \in \lim_P T$.
Let $x \in P$ and note that $s_x \in T_x=\lim_{(0_{\widehat P},x]} T$.
Since $T_{x}$ is an ASL on $(0_{\widehat P},x]^{\op}$ there exist
unique scalars $\lambda^x_{m,s} \in K$ such that
$$
s_x
=
\sum_{\Mcc \in \N^{(0_{\widehat P},x]} \cap \stand} \lambda^x_{\Mcc,s}
  \phi_{x}(\Mcc).
$$
For every chain $\Mcc$ in $P$ we set $x_{\Mcc} = \max(\supp(\Mcc))$.
Note that by the assumptions we have that
$\lambda^x_{\Mcc,r} = \lambda^y_{\Mcc,r}$
for all $x,y \in P$ such that $y\geq x$
and $x_\Mcc  \leq x,y$.
We define $r_{\Mcc,s}= \lambda^{x_{\Mcc}}_{\Mcc,s}$.
Let
$$
f_s=
\sum_{\Mcc \in \stand}
r_{\Mcc,s} \Mcc \in K[P].
$$
Then by the choice of the coefficients $r_{\Mcc,s}$
$$
\phi(f_s) = s \in \lim_P T.
$$
This shows the first claim.
Observe that if $x,y \in P$ are incomparable,
then
$x_\Mcc>x,y$
for all $\Mcc$ such that $r_{\Mcc,\phi(xy)}\neq 0$. Thus ASL2 follows
ones
we have proved
that the standard monomials are $K$-linearly independent in $\lim_P T$.
Assume that
$$
0=
\sum_{\Mcc \in \stand}
a_{\Mcc} \phi(\Mcc) \in \lim_P T.
$$
Choose ${\Ncc} \in \stand$ and
consider the projection $\pi(T)_{x_\Ncc}\colon \lim_P T \to
T_{x_{\Ncc}}$.
Hence
$$
0
=
\sum_{\Mcc \in \stand}
 \pi(T)_{x_\Ncc}( a_{\Mcc} \phi(\Mcc) )
=
\sum_{\Mcc \in \N^{(0_{\widehat P}, x_{\Ncc}] } \cap \stand}
a_{\Mcc} \phi_{x_\Ncc}(\Mcc) \in T_{x_\Ncc}.
$$
Since $T_{x_\Ncc}$ is an ASL $(0_{\widehat P},x_\Ncc]^{\op}$, we get that $a_{\Ncc}=0$.
This shows that the standard monomials are indeed
$K$-linearly independent in $\lim_P T$ and thus we have proved
that ASL1 holds. As noted above also ASL2 holds and this shows
that $\lim_P T$ is an ASL on $P^{op}$.
\end{proof}

\begin{ex}[Locally distributive lattices]
\label{exsimplicial}
We refer to \cite{ST99} for the notion of distributive lattice.
A finite poset $P$ is called
a {\em locally distributive lattice},
if $P$ has a terminal element $0_P$ and
for all elements $z \in P$ the interval
$[0_P,z]$ is a distributive lattice.
For example a
simplicial posets
as considered in \cite{DU97}, \cite{MA04} and \cite{ST96},
i.e. $P$ has $0_P$ and
for all $z \in P$ the interval
$[0_P,z]$ is a boolean algebra,
are locally distributive lattices.
Let $z \in P$.
Since $[0_P,z]$
is a distributive lattice
it follows from a result of Hibi (see \cite{HI87})
that
$$
T_z=K[[0_P,z]]/(x\cdot y-x\wedge y \cdot x\vee y)
$$
is a graded ASL on $P^{\op}$ with straightening relation
$r_z(x,y)=x\wedge y \cdot x\vee y$ (which is an integral domain).
Here $\wedge$ and $\vee$ are the meet and join in the distributive lattice
$[0_P,z]$.
Note that $x\wedge y$ is the meet of $x$ and $y$ in $P$ and $x\vee
y$ is a minimal upper bound of $x,y $ in $P$ depending on $z$.

Let $R=K[P]$.
The restriction homomorphisms
$T_{z'z}: T_z \to T_{z'}$ for $z,z' \in P$
with $z'\leq z$,
define a $\Z$-graded $KP$-algebra.
We call $\lim T$
the {\em generalized Hibi ring} associated to $P$.
By Theorem \ref{aslconstr} $\lim_P T$ is an ASL on $P^{\op}$
and it is easy to see that for incomparable $x,y \in P$ the
straightening
relations are
\begin{eqnarray*}
r(x,y) = x\wedge y \sum_{z } z
\end{eqnarray*}
where $z$ ranges over all minimal upper bounds for $x,y$.
The sum is defined to be zero if there are no such elements.
\end{ex}

Note that if $R$ is a {\em graded ASL}, i.e.
$R$ is a graded $K$-algebra
and all elements of $P$ are homogeneous and have positive degree,
then it is known that
the straightening relation (\ref{asleqn})
give a presentation of $R$ (see Proposition 1.1 in \cite{COEICL}).
Observe that this is not true in the general case
(see \cite{MI93} and \cite{TR91} for counterexamples).

Together with this observation we proved in
Example \ref{exsimplicial}
in fact the following corollary:

\begin{cor}
\label{aslring}
Let $P$ be a locally distributive lattice and
let $T$ be the $KP$-algebra constructed in
\ref{exsimplicial}.
Then
$\lim_P T$ is a $\Z$-graded ASL on $P^{\op}$
and for $\lim_P T=K[P]/I_{T}$
we have that
the ideal
$I_{T}$
is generated
by
\begin{eqnarray*}
x y &-& x\wedge y \sum_{z } z
\end{eqnarray*}
where $z$ ranges over all minimal upper bounds for $x,y$.
The sum is defined to be zero if there are no such elements.
\end{cor}
Hence we obtain exactly the poset ring $\Tilde{A}_P$ constructed by
Stanley in \cite{ST91} for simplicial posets with the same ASL structure.
One immediately gets the following result:

\begin{cor}
Let $P$ and $T$ be as in Corollary \ref{aslring}.
If $P$ is a Cohen--Macaulay poset, i.e.\ $\Delta(P)$ is Cohen--Macaulay simplicial complex,
then
$\lim_P T$ is a Cohen--Macaulay ring.
\end{cor}
\begin{proof}
Since $\lim_P T$ is $\Z$-graded and
the elements of $P$
have positive degree,
it follows from
\cite[Cor. 7.2]{COEICL}
that
$\lim_P T$ is a
Cohen--Macaulay ring.
\end{proof}


\section{Quotients of rings of global sections}
\label{quotients}
Let
$P$ be a finite poset and
$T$ be a flasque $RP$-algebra.
Given an open subset $Q$ of $P$,
the restriction of $T$
to $Q$ induces a surjective homomorphism
$$
 \lim_P T \to \lim_Q T.
$$
Thus $\lim_Q T$ is a quotient of $\lim_P T$
and we are interested in the relationship between
ring properties of these two rings. In general there is no close
connection,
but in special situations there is more hope.

The goal of this section is to
generalize so-called
{\em rank--selection theorems} in the tradition of
Duval \cite{DU97},
Hibi \cite{HI91}, Munkres \cite{MU84}
and Stanley \cite{ST79}
to our situation. We follow the ideas of Hibi \cite{HI91}
and at first we define:

\begin{defn}
\
\begin{enumerate}
\item
Let $X$ be a subset of $P$. We denote by $P_X$ the sub-poset
$P\setminus \bigl(\bigcup_{x \in X} [x,1_{\widehat P})\bigr)$.
\item
For $x \in P$ the {\em star} of $x$ in $P$ is the sub-poset
$\str_P(x)$ consisting of the elements $y \in P$ such that
there exists $z \in P$ with $x, y\leq z$.
\item
A non-empty subset $X \subseteq P$
is called {\em excellent} if for every maximal element $z$ of $P$
there exists at most one $x \in X$ with $x \le z$.
\end{enumerate}
\end{defn}

Given an element $x \in P$
we let $I_x$ denote
the kernel of the homomorphism
$$
\lim_P T \to \lim_{P\setminus [x,\widehat 1)} T .
$$

\begin{lem}
\label{shorthelper}
Let $R$ be a $\Z^m$-graded ring,
$X$ be an excellent subset of a finite poset $P$ and
$T$ be a flasque $\Z^m$-graded $RP$-algebra.
Then
there exists a short exact sequence
\begin{displaymath}
0
\to
\bigoplus_{x \in X} I_x
\to
\lim_P T
\to
\lim_{P_X} T
\to
0
\end{displaymath}
of $\Z^m$-graded $R$-modules.
\end{lem}
\begin{proof}
Let $x \in X$.
The homomorphism
$$
\bigoplus_{x \in X} I_x
\to
\lim_P T
$$
is the sum of the inclusions $I_x \to \lim_P T$.
This homomorphism is injective
because $X$ is excellent.
It follows from the definition of $P_X$ that
the sequence in question is a zero-sequence.
It is right-exact,  since $T$ is flasque.
We leave it as an exercise to the reader to check middle-exactness.
\end{proof}
\begin{defn}
A poset $P$ is {\em locally graded}
if the maximal chains of $(0_{\widehat P},x]$ are
of equal cardinality for every $x \in P$ .
\end{defn}
We can now apply the main result in \cite{BBR05}.
\begin{lem}
\label{CMhelper}
Let $(R,\mm)$ be a $\Z^m$-graded local ring with $K = R_0$ a field and
let $P$ be a finite locally
graded poset.
Let
$T$ be a flasque $\Z^m$-graded $RP$-algebra such that the homomorphism
$R \to \lim_P T$ is surjective.
Assume that:
\begin{enumerate}
\item
$T_y$ is Cohen--Macaulay of dimension $\rk(y)$ for every $y \in P$.
\item
For every $y \in P$
we have that
$$
\widetilde H^{i}((y,1_{\widehat P});K)=0 \text{ for }i \neq \rk((y,1_{\widehat P})).
$$
\end{enumerate}
Then $\lim_P T$ is Cohen--Macaulay with $\dim \lim_P T = \rank(P)$.
In particular, if $P$ is a Cohen--Macaulay poset then
the assumption (ii) is satisfied.
\end{lem}
\begin{proof}
By \cite[Theorem 6.2]{BBR05}
we have that
$$
H^i_\mm(\lim_P T)
\cong
\bigoplus_{y \in P}
\widetilde H^{i - \rk(y) -1}((y,1_{\widehat P});K)
\tensor_K
H^{\rk(y)}_\mm(T_y)
$$
as $\Z^m$-graded $K$-vector spaces.
By assumption (ii)
we have that
$$
\widetilde H^{i - \rk(y) -1}((y,1_{\widehat P});K)
=0
$$
$$
\text{ for }
i - \rk(y) -1 \neq \rk((y,1_{\widehat P}))= \rk(P)-\rk(y)-1.
$$
The latter condition is equivalent to  $i \neq \rk(P)$.
Hence
$$
H^i_\mm(\lim_P T)=0
\text{ for }
i \neq \rk(P)
$$
and thus $\lim_P T$ is Cohen--Macaulay of dimension $\rk(P)$.
\end{proof}

We are ready to present the main result of this section which
generalizes a result of Hibi \cite[Theorem 1.7]{HI91}.

\begin{thm}
\label{remCMthm}
Let $(R,\mm)$ be a $\Z^m$-graded local ring with $K = R_0$ a field and
let $X$ be an excellent
subset
of a
finite locally
graded poset $P$.
Let
$T$ be a flasque $\Z^m$-graded $RP$-algebra
with $\dim \lim_P T = \rank(P)$, and such that the homomorphism $R \to
\lim_P T$ is surjective.
Assume that:
\begin{enumerate}
\item
$T_y$ is Cohen--Macaulay of dimension $\rk(y)$ for every $y \in P$.
\item
For every $x\in X$ we have that $\rank(\str(x)) = \rank(P)$ and
$\rank(\str(x) \cap P_X) = \rank(P) -1$.
\item
For every $x \in X$ and $y \in \str(x)$ we
that
$\widetilde H^{i}((y,1_{\widehat \str(x)});K)=0$
for $i \neq \rk((y,1_{\widehat \str(x)}))$.
\item
For every $x \in X$ and  $y \in \str(x)\cap P_X$
we have
that
$\widetilde H^{i}((y,1_{\widehat{\str(x)\cap P_X }});K)=0$
for $i \neq \rk((y,1_{\widehat{\str(x)\cap P_X }}))$.
\end{enumerate}
Then the natural homomorphism
$$
H^i_{\mm}( \lim_P T ) \to H^i_{\mm}(\lim_{P_X} T )
$$
is surjective for $i = \rank(P)$,
it is injective for $i = \rank(P)-1$, and otherwise it is an isomorphism.
\end{thm}
\begin{proof}
For every $x \in X$ the sets
$\str(x)$ and $\str(x)\cap P_X$ are open subsets of $P$
and thus they are locally graded posets and the ranks of elements in these
poset coincide with the ones in $P$.
Let $J_x$ be the kernel of
the surjective homomorphism
$\lim_{\str(x)} T \to \lim_{\str(x)\cap P_X} T$.
Since
$\str(x)\cap P_X=\str(x)\setminus [x,1_{\widehat P})$
we get a commutative diagram of $R$-modules of the form
\begin{center}
$$
\begin{CD}
0 @>>> I_x @>>> \lim_P T @>>> \lim_{P\setminus [x,1_{\widehat P})} T @>>> 0\\
  & & @VVV  @VVV  @VVV\\
0 @>>> J_x @>>> \lim_{\str(x)} T @>>> \lim_{\str(x)\cap P_X} T @>>> 0.
\end{CD}
$$
\end{center}
Note that $\lim_P T \to \lim_{\str(x)}T$ and
$\lim_{P\setminus [x,1_{\widehat P})} T \to \lim_{\str(x)\cap P_X} T$
are surjective and
that $I_x \to J_x$ is an isomorphism.

By Lemma \ref{CMhelper} the rings $\lim_{\str(x)}T $
and
$\lim_{\str(x)\cap P_X} T $
are Cohen--Macaulay of dimension
$\rk(\str(x))=\rk(P)$ and
$\rk(\str(x)\cap P_X)=\rk(P)-1$ respectively.
The long exact local cohomology sequence induced by the second row above
shows that $J_x$ and thus $I_x$ are Cohen--Macaulay of dimension
$\rk(P)$.

Next we consider the short exact sequence
\begin{displaymath}
0
\to
\bigoplus_{x \in X} I_x
\to
\lim_P T
\to
 \lim_{P_X} T
\to
0
\end{displaymath}
and read the statement of the theorem off from the associated long
exact local cohomology sequence.
\end{proof}
If $X$ is an excellent subset of a finite locally graded poset $P$,
then $P_X$ is
called {\em hereditary}
if
the assumptions (ii)--(iv) of Theorem \ref{remCMthm}
are satisfied.
In view of the above result and a theorem of Hochster implying that
the normal monoid ring $K[D\cap \Z^d]$ is Cohen--Macaulay of
dimension $d$ for every $d$-dimensional
rational pointed cone $D$ in $\R^d$,
we obtain the following

\begin{cor}
\label{coni}
Let $P(\Sigma)$ be the face poset of a rational pointed $d$-dimensional fan $\Sigma$
in $\R^m$, let $R = \lim_P T$
and let $T(\Sigma)$ be the associated $KP(\Sigma)$-algebra of Example \ref{exfan}.
If $X$ is a hereditary subset of $P(\Sigma)$
then
the homomorphism
$$H^i_{\mm}(\lim_{P(\Sigma)} T) \to
H^i_{\mm}(\lim_{P(\Sigma)_X} T)$$
is surjective for $i = d$, it is injective for $i = d-1$, and otherwise it is an isomorphism.
\end{cor}

A sub-poset $Q$ of $P$ is called {\em $n$-hereditary}
(or hereditary for short)
if there exists a sequence
$X_1, X_2,\dots,X_n$ of
subsets of $P$ and a sequence
$Q = P_1 \subseteq P_2 \subseteq \dots
\subseteq P_{n+1} = P$ of sub-posets of $P$ such that $X_i$ is an excellent subset of
$P_{i+1}$ and $P_i = (P_{i+1})_{X_i}$ hereditary in $P_{i+1}$ for
$i = 1,\dots,n$.

\begin{cor}
\label{nicetohave}
Let $Q$ be an $n$-hereditary sub-poset of a locally graded
poset $P$.
Let $(R,\mm)$ be a $\Z^m$-graded local ring with $R_0 =K$ a field and
let
$T$ be a flasque $\Z^m$-graded $RP$-algebra
such that
the homomorphism $R \to \lim_P T$ is surjective and such that
$T_y$ is Cohen--Macaulay of dimension $\rk(y)$ for every $y \in P$.
Then
$$
H^i_\mm( \lim_P T)
\cong
H^i_\mm(\lim_Q T)
\text{ for }
i=0,\dots, \rk(P)-1-n
$$
as $\Z^m$-graded $\lim_P T$-modules.
In particular,
if $\lim_P T$ is Cohen--Macaulay and $\rk(Q) = \rk(P)-n$,
then $\lim_Q T$ is Cohen--Macaulay.
\end{cor}

The {\em $i$th skeleton} of a poset $P$ is
the sub-poset
$$
P^{(i)}=\{y \in P:\rk(y) \leq i\}.
$$

\begin{ex}
\label{fanskeleton}
Let $P=P(\Sigma)$ be the face poset of a rational pointed fan $\Sigma$
in $\R^m$ and let
$T=T(\Sigma)$ be the $KP(\Sigma)$-algebra of Example
\ref{exfan}.
We claim that
$P^{(i)}$ is a hereditary sub-poset of $P^{(i+1)}$
with respect to the excellent subset
$$
X^{(i)}=\{C \in P^{(i+1)}: \rk(C) =i+1 \}.
$$
For
$C \in X^{(i)}$
we have that $\str(C) \cap P^{(i+1)}$ contains an initial element,
and thus it is acyclic.
For every $C \in P$ the poset
$\str(C)\setminus \{C\}$
is the face poset of the boundary of a cone and thus
Cohen--Macaulay. This shows the claim.

We obtain immediately
from Corollary \ref{nicetohave}
that
$$
\depth \lim_P T
=
\max\{i \in \N : i \le \rank(P) \text{ and } \lim_{P^{(i)}} T \text{
  is Cohen--Macaulay} \}.
$$
If $P$ is the face poset of a simplicial complex $\Delta$ and $T$ is
the Stanley--Reisner $KP$-algebra associated to $\Delta$, then this
is a well-known result.
\end{ex}

\section{Presentations of sectioned $KP$-algebras}
\label{sectioned}
Recall that for a set $F$ we denote by $K[F]$
the polynomial ring on the variables $a \in F$ and that the {\em
  support} of a monomial $\Mcc \in K[F]$ is the set
$\supp(\Mcc)=\{a \in F : \Mcc(a)>0\}$.
Given an inclusion $F_x \subseteq F$ of finite sets we denote by
$\pi_x
\colon K[F] \to K[F_x]$ the natural projection with $\pi_x(a) = a$ if
$a
\in F_x$ and $\pi_x(a) = 0$ otherwise. There is a $KP$-algebra
$T_{\Pi}$ associated to every system $\Pi = (F_x)_{x \in P}$ of
subsets of a finite set
$F$ with $F_x \cap F_y = \cup_{z \le x,y} F_z$ for every $x,y \in P$
defined by $({T^{\Pi}})_x = K[F_x]$ and
$({T^{\Pi}})_{xy}(f) = \pi_x(f)$ for $x \le y$ in $P$ and $f \in K[F_y]
\subseteq K[F]$.
The projections $\pi_x \colon K[F] \to K[F_x]$ induce a surjective
homomorphism $\pi \colon K[F] \to \lim_P T^{\Pi}$
with kernel $\Ker(\pi) = \bigcap_{x \in P} \Ker(\pi_x)$ equal to the
ideal
generated by the set of monomials $\Mcc$ in $K[F]$ such that $\supp(\Mcc)$ is
not contained in any of the sets $F_x$.
Thus the ring $\lim_P T^{\Pi}$ is isomorphic to the Stanley--Reisner
ring of
the simplicial complex $\Delta_{\Pi}$ consisting of the subsets
$G$ of $F$ contained in $F_x$ for some $x \in P$.

\begin{defn}
Let $T$ be a $KP$-algebra.
We call a pair $(\Pi,p)$ of a system $\Pi = (F_x)_{x \in P}$ of
subsets of a finite set
$F$ with $F_x \cap F_y = \cup_{z \le x,y} F_z$ for every $x,y \in P$
and a surjection $p \colon {T^{\Pi}} \to T$ of $KP$-algebras
a {\em presentation} of $T$. The homomorphisms induced by $p$ are
denoted $p_Q \colon \lim_Q
T^{\Pi} \to \lim_Q T$ for $Q$ open in $P$ and $p_x \colon (T^{\Pi})_x
\to T_x$ for $x \in P$.
\end{defn}
Given a $KP$-algebra $T$ and $x \in P$, the
structure-homomorphism from
$\lim_P T$ to $T_x$ is denoted $\pi(T)_x \colon \lim_P T \to T_x$.
\begin{defn}
A {\em sectioned $KP$-algebra} $(T,\iota(T))$ consists of a
$KP$-algebra $T$ and
$K$-linear homomorphisms $\iota(T)_x \colon T_x \to \lim_P T$ with
$\pi(T)_x \circ \iota(T)_x = \id_{T_x}$ for  $x \in P$.
\end{defn}
\begin{ex}[Monomial ideals]
  Let $I$ be a monomial ideal in the polynomial ring $S =
  K[x_1,\dots,x_n]$ with irredundant irreducible decomposition $I =
  \mm^{b^1} \cap \dots \cap \mm^{b^t}$. We consider the $KP$-algebra
  $T = T(I)$ on the poset $P = P(I)$ of Example \ref{exmulti}. For $D
  \subseteq C \subseteq [t]$, there is a preferred section
  $S/\sum_{i \in D} \mm^{b^i} \to S/\sum_{i \in C} \mm^{b^i}$ of the projection
  $T_{DC}$ induced by a $K$-linear homomorphism $S \to S$ acting
  either as the identity or as zero on a monomial of $S$. The
  $KP$-algebra $T$ is sectioned with $\iota(T)_C$ such that the
  composition $\pi(T)_D \circ \iota(T)_C$ is given by the composition
  \begin{displaymath}
    T_C = S/\sum_{i \in C} \mm^{b^i} \to S/\sum_{i \in C \cap D}
    \mm^{b^i} \to S/\sum_{i \in D} \mm^{b^i},
  \end{displaymath}
  where the first homomorphism is the section described above and the second
  homomorphism is the natural projection.
\end{ex}

The $KP$-algebras of the form $T^{\Pi}$ considered above are sectioned
$KP$-algebras with
$K$-linear sections $\iota(T^{\Pi})_x \colon {T^{\Pi}}_x \to \lim_P
{T^{\Pi}}$
defined by $\iota(T^{\Pi})_x(f) = (p_y(\iota_x(f)))_{y \in P}$, where
$\iota_x \colon K[F_x] \to K[F]$ denotes the natural inclusion.
\begin{defn}
We call $(\Pi,p)$ a {\em presentation} of the sectioned $KP$-algebra
$T$ if it is a presentation of $T$ considered as a $KP$-algebra and
the identity
$p_P \circ \iota(T^{\Pi})_x =
\iota(T)_x \circ p_x$ holds for every $x \in P$.
\end{defn}
The notation of the above definition can be summarized in the
following commutative diagram:
\begin{displaymath}
  \begin{CD}
   K[F] @= K[F] @= K[F] \\
   @A{\iota_x}AA @V{\pi}VV @V{\pi_x}VV \\
   K[F_x] @>{\iota(T^{\Pi})_x}>> \lim_P T^{\Pi}
   @>{\pi(T^{\Pi})_x}>> K[F_x] \\
   @V{p_x}VV @V{p_P}VV @V{p_x}VV \\
    T_x @>{\iota(T)_x}>> \lim_P T @>{\pi(T)_x}>> T_x.
  \end{CD}
\end{displaymath}
Let us record some facts.
\begin{lem}
If $(\Pi,p)$ is a presentation of a sectioned $KP$-algebra $T$, then
\begin{enumerate}
\item
$
0 =
\bigcap_{x \in P} \Ker(\pi(T)_x)
\subset \lim T.
$
\item
$
\Ker(p_P) = \bigcap_{x \in P} \Ker(\pi(T)_x \circ p_P).
$
\item
$
\Ker(\pi(T)_x \circ p_P \circ \pi) =
\Ker(p_x \circ \pi_x) = \iota_x(\Ker(p_x)) + (F
\setminus F_x).
$
\item
$
\iota_x(\Ker(p_x)) =
\Ker(p_x \circ \pi_x) \cap \iota_x(K[F_x]) = \Ker(\pi(T)_x \circ p_P
\circ \pi)
\cap
\iota_x(K[F_x]).
$
\end{enumerate}
\end{lem}

The following theorem is our main result in this section.
\begin{thm}
\label{presentation}
Let
$(\Pi,p)$ be
a presentation of a sectioned $KP$-algebra $T$.
We have that
\begin{enumerate}
\item
$\Ker(\pi) = \bigcap_{x \in P} \Ker(\pi_x)=I_{\Delta_\Pi}$ is the
  Stanley--Reisner
ideal of $\Delta_\Pi$ in $K[F]$.
\item
$
\Ker(p_P \circ \pi)=
\sum_{x \in P} (\iota_x (\Ker(p_x)))
+
\bigcap_{x \in P} \Ker(\pi_{x}).
$
\end{enumerate}
In particular,
$$
\lim_P T \cong
K[F]/
\sum_{x \in P} I_x + I_{\Delta_\Pi},
$$
where
$I_x=(\Ker(p_x \circ \pi_{x}) \cap \iota_x(K[F_x]))$ for $x \in P$.
\end{thm}
\begin{proof}
Since we have explained (i) above we only need to prove (ii).
Clearly the right hand side is contained in $\Ker(p_P \circ \pi)$.
  Given $f \in \Ker(p_P \circ \pi)$ we consider the set $P(f)$
  consisting of the elements $x
  \in P$ such that $f$ contains a monomial $\Mcc$ with $\supp(\Mcc)
  \subseteq F_x$.
  We prove the
  opposite inclusion by induction on the cardinality of the set
  $P(f)$.
  First note that $P(f) = \emptyset$ implies $f \in \bigcap_{x \in P}
    \Ker(\pi_x)$.
  If $P(f) \neq \emptyset$ we choose an element $x
  \in P(f)$. Since $\iota_x \pi_x(f) \in \iota_x(\Ker(p_x))$
  the
  element $f' = f - \iota_x \pi_x(f)$ is an element of
  $\Ker(p_P \circ \pi)$ with $P(f')$ a proper subset of $P(f)$. This takes
  care of the induction step.
\end{proof}

In some cases
there is a  relation between
the simplicial complex $\Delta_{\Pi}$ and the
$K$-algebra $\lim_P T$.
In fact, we have:

\begin{cor}
\label{presnice}
Let $(\Pi,p)$ be a presentation of a sectioned flasque $KP$-algebra
$T$ consisting of
$\Z^m$-graded Cohen--Macaulay rings such that
$d_x<d_y$ and $f_x<f_y$ for $x < y$ in  $P$,
where
$d_x=\dim(T_x)$ and
$f_x=|F_x|$.
We have that
\begin{eqnarray*}
\depth \lim_P T
&=&
\min\{ i \in \N :
H^{i+(f_x-d_x)}_\mm(K[\Delta_\Pi])\neq 0
\text{ for some } x \in P
\},
\\
\dim \lim_P T
&=&
\max\{ i \in \N :
H^{i+(f_x-d_x)}_\mm(K[\Delta_\Pi])\neq 0
\text{ for some } x \in P
\}.
\end{eqnarray*}
In particular, if there exists an $n \in \Z$
such that
$\dim (K[F_x]) = \dim(T_x)+n$ for every $x \in P$ and such that
$\dim(K[\Delta_{\Pi}]) = \dim(\lim_P T) + n$, then
$\lim_P T$ is Cohen--Macaulay if and only if
$\Delta_{\Pi}$ is Cohen--Macaulay.
\end{cor}
\begin{proof}
The $K$-algebra $\lim_P T^{\Pi}$ is isomorphic to the
Stanley--Reisner ring of $\Delta_{\Pi}$.
The statements are a direct consequence of
\cite[Corollary 6.3]{BBR05}.
\end{proof}

\begin{ex}[Toric face rings]
Let $\Sigma$ be a rational pointed fan in $\R^m$ and choose affine
monoids $M_C$ for $C \in \Sigma$ as in Example \ref{exfan}.
The $KP$-algebra $T$ of Example
\ref{exfan} with $T_C = K[M_C]$ is a sectioned $KP$-algebra with
$\iota(T)_C \colon T_C \to \lim_P T$ determined by requiring the
composition $\pi(T)_D \circ \iota(T)_C$ to be the composition of the
inclusion $K[M_C] \subseteq K[M_E]$ and the face projection $K[M_E]
\to K[M_D]$ if $C$ and $D$ are faces of a common cone $E$ in $\Sigma$, and
requiring $\pi(T)_D \circ \iota(T)_C$ to be zero otherwise.
Denoting the Hilbert basis of $M_C$ by $F_C$, the surjections
$p_C \colon K[F_C] \to K[M_C]$ define a
presentation $(\Pi,p)$ of the sectioned $KP$-algebra $T$ with
$\Pi = (F_C)_{C \in \Sigma}$.
In the particular case where $M_C = \Z^m \cap C$ for every $C \in
\Sigma$, the $KP$-algebra $T$ is denoted $T(\Sigma)$, and
$T(\Sigma)_C$ is Cohen--Macaulay of dimension $\dim(C)$ since it is
the monoid ring of a normal affine monoid.
Note that Corollary \ref{presnice} applies to
$T(\Sigma)$.
\end{ex}

Finally, we want to compute the initial ideals of
the presentation ideal of $\lim_P T$ with respect to
weight orders $\ini_\omega$ on $K[F]$
induced by a map  $\omega \colon F \to \R$.

\begin{thm}
Let $(\Pi,p)$ be a presentation of a sectioned $KP$-algebra $T$.
For a map $\omega \colon F \to \R$
we have that
$$
K[F]/{\ini_{\omega}(\Ker(p_P \circ \pi))} \cong \lim_P
K[F_x]/\ini_{\omega}(\Ker(p_x)).
$$
\end{thm}
\begin{proof}
  The $KP$-algebra $\ini_{\omega} T$ with $\ini_{\omega} T_x =
  K[F_x]/\ini_{\omega}(\Ker(p_x))$ is presented by
  $(\Pi,q)$, where
  $q_x$ is the projection $K[F_x] \to K[F_x]/\ini_{\omega}(\Ker(p_x
  ))$.
  Hence by Theorem \ref{presentation}
  \begin{displaymath}
    \Ker(q_P \circ \pi) = \left(\sum_{x \in P} \iota_x(\Ker(q_x )) \right) +
    \bigcap_{x \in P} \Ker(\pi_x).
  \end{displaymath}
  Clearly $\iota_x(\Ker(q_x )) =
  \iota_x(\ini_{\omega}(\Ker(p_x )))$ is contained
  in $\ini_{\omega}(\Ker(p_P \circ \pi))$, and thus the identity on
  $K[F]$ induces
  a homomorphism $K[F]/\Ker(q_P \circ \pi) \to
  K[F]/\ini_{\omega}(\Ker(p_P \circ \pi))$.
  On the other hand, $\pi_x(\ini_{\omega}(\Ker(p_P \circ \pi))) \subseteq
  \ini_{\omega}(\Ker(p_x ))$, and thus the projection $K[F] \to
  K[F_x]/\ini_{\omega}(\Ker(p_x ))$ induces a homomorphism
  $$K[F]/\ini_{\omega}(\Ker(p_P \circ \pi)) \to
  K[F_x]/\ini_{\omega}(\Ker(p_x )).$$
  These homomorphisms in turn
  assemble to a homomorphism
  \begin{displaymath}
    K[F]/{\ini_{\omega}(\Ker(p_P \circ \pi))} \to \lim_P
    K[F_x]/\ini_{\omega}(\Ker(p_x ))
    \cong K[F]/\Ker(q_P \circ \pi).
  \end{displaymath}
We leave it as an exercise to the reader that this is an
inverse isomorphisms to $K[F]/\Ker(q_P \circ \pi) \to
K[F]/\ini_{\omega}(\Ker(p_P \circ \pi))$.
\end{proof}

\end{document}